\documentclass[opre,nonblindrev]{informs3_no_remark} 

\OneAndAHalfSpacedXI 

\usepackage{natbib}
 \bibpunct[, ]{(}{)}{,}{a}{}{,}%
 \def\bibfont{\small}%

 \usepackage{mathtools}
 \usepackage{tikz,calc}
 \usetikzlibrary{decorations.markings}
 \usepackage{soul}
 \setul{}{0.75pt}
 \setuloverlap{1.5pt}

\usepackage{subfigure}
\usepackage{mwe}
\usepackage{lscape,booktabs,longtable}
\usepackage{amsmath}
\usepackage[pdftex,colorlinks=true,urlcolor=blue,citecolor=blue,pdfstartview=FitH]{hyperref}
\usepackage{comment}
\usepackage{enumerate}
\usepackage[ruled,vlined]{algorithm2e}
\usepackage[graphicx]{realboxes}

\newcommand{\E}[1]{{\mathrm{E}\left[#1\right]}}

\newcommand{\PP}[1]{\mathrm{P}\left( #1 \right)}

\TheoremsNumberedThrough     
\ECRepeatTheorems

\EquationsNumberedThrough    

\begin{document}


\RUNAUTHOR{Daw and Goel}

\RUNTITLE{On the Duration of Markovian Hawkes Clusters}

\TITLE{\textsf{On the Cluster Duration Distribution of the Markovian Hawkes Process}
}
\ARTICLEAUTHORS{%
\AUTHOR{Andrew Daw and Amer Goel}
\AFF{Dept.~of Data Sciences and Operations, University of Southern California Marshall School of Business, Los Angeles, CA} 

} 


\ABSTRACT{
The duration of a Hawkes process cluster, meaning the time from first point to last point in an entirely endogenously generated progeny stream, has been stubbornly resistant to analysis; even the mean has not been known in closed form.
Recently, uniformly random parking functions have been found to be hidden spines of these history-driven clusters, offering a probabilistic combinatorial structure within self-excitement that has opened new avenues for analysis. Here, we leverage and refine this decomposition to provide a closed-form expression for the mean cluster duration in the case of an exponential excitation kernel. Specifically, we reduce from uniformly random parking functions to non-uniform preference counts. This change of probabilistic combinatorial structure in the decomposition leads us to both a Borel-based identity and an explicit integral-based expression for the mean duration.
}
\KEYWORDS{Hawkes process, cluster duration, end-state conditioning, generalized parking functions, preference counts} 

\maketitle


%


\vspace{-0.12in}

\section{Introduction}
\label{sec: intro}


Originally introduced in \citet{hawkes1971spectra}, the Hawkes point process is a \emph{self-exciting} model, in which the history of points so far is an endogenous driver of future points to come. This history-driven process has enjoyed myriad applications over the past half-century, ranging from seismology to neurology, from financial trading to risk, and from pricing to public health \citep[e.g.,][]{bacry2015hawkes,brown2002time,liao2026tides,ogata1988statistical,xu2014path,yuan2026dynamic}. Much of the fundamental understanding of this model comes from the seminal \citet{hawkes1974cluster} decomposition, which distills the model into two separate streams: exogenous Poisson arrivals and exclusively endogenous clusters that descend as self-excited progeny from the exogenous points. Hence,  we can isolate the self-excitement to a purely history-driven \emph{Hawkes cluster model} which has the stochastic intensity and counting process pair $(\lambda_t, N_t)$ governed by
\begin{align}
\lambda_t
&=
\int_0^t g(t-u)\mathrm{d}N_u
=
\sum_{i=0}^{N_t-1}
g(t-\tau_i)
,
\label{intensityDefEq}
\end{align}
where $g : \mathbb{R}_+ \to \mathbb{R}_+$ is the \emph{excitation kernel} and $0 = \tau_0 < \tau_1 < \tau_2 < \dots$ are the \emph{cluster points}. If we specify $g(\cdot)$ to be $g(t) = \rho e^{-t}$ for some $\rho \in (0,1)$, then we receive the cluster model for the \textit{Markovian Hawkes process}, which is likely the most widely used and well-studied form of self-excitement \citep{oakes1975markovian}. 

One of the reasons that the \citet{hawkes1974cluster} decomposition has become such a foundational tool for studying self-excitement is that time-agnostic quantities become readily available. By connection to branching processes with Poisson-distributed offspring, a great deal can be known about the model at the count-level: given $\rho < 1$, then the cluster will have only finitely many points almost surely ($\PP{N < \infty} = 1$ for $N := \lim_{t \to \infty} N_t$); moreover, the distribution of the \emph{size} of the cluster is known in closed form ($N \sim \mathsf{Borel}(\rho)$). 
However, for a model that is inarguably temporal, size only says so much. The corresponding first-order temporal quantity, the \emph{duration} of the cluster $\tau := \tau_{N-1}$, is undoubtedly at least as important in practice as the size, but, unfortunately, chronological understanding of Hawkes clusters has been much more elusive. As said by \citet[][pg.~631]{moller2005perfect}, the distribution of the cluster duration ``is unknown even for the simplest examples of Hawkes processes.''


The analytical elusiveness of the Hawkes cluster duration is not from a lack of effort. For instance, \citet{moller2005perfect} themselves provided an integral equation for $F(t) = \PP{\tau \leq t}$ \citep[Prop. 2 therein, which the authors describe as a generalization of Theorem 5 from][]{hawkes1974cluster}. Sticking with the simple case that $g(t) = \rho e^{-t}$, this integral equation can be seen to be equivalent to the ordinary differential equation (ODE) $F'(t) = F(t)\left(\rho F(t) - \rho - \log(F(t))\right)$ with initial condition $F(0) = \PP{N = 1} = e^{-\rho}$; as noted in Example 4 of \citet{moller2005perfect}, this ODE is not analytically solvable. Building on \citet{moller2005perfect}, \citet{graham2021regenerative} showed that $\E{\tau} \leq 1/(1-\rho)$ and that, moreover, $\tau$ is stochastically dominated by an $\mathsf{Exp}(1-\rho)$ random variable. As another example of the techniques  derived from \citet{hawkes1974cluster}, \citet{chen2021perfect} constructed an almost sure bound on $\tau$ through the ``total birth time,'' which, for $g(t) = \rho e^{-t}$, is equivalent to the sum of a $\mathsf{Borel}(\rho)$-distributed number of mutually independent $\mathsf{Exp}(1)$-distributed random variables. (This random sum yields a similar mean bound, $\E{\tau} \leq \rho/(1-\rho)$.) 
These examples are meant to be representative but not exhaustive, and additional bounds on the mean or the tail of the duration have been provided by works such as \citet{bremaud2002rate,chen2020perfect,costa2020renewal,reynaud2007some}. 
Recently, and quite relevantly for this paper's pursuit, there has been progress on understanding the conditional distribution of $\tau$ given $N$. Via an extension of the random time change theorem, \citet{daw2023conditional} showed that, for $g(t) = \rho e^{-t}$, $\tau$ is equivalent in distribution to a sum of conditionally independent (but not identical) exponential random variables with rates determined by certain probabilistic combinatorial objects, and these conditional connections were subsequently refined in \citet{daw2026asymmetries}.
Nevertheless, to the best of our knowledge, understanding of the unconditional distribution remains lacking: even the simplest performance metric, $\E{\tau}$, is not known in closed form for even the simplest excitation kernel, $g(t) = \rho e^{-t}$.


\textbf{In this paper,  we explicitly characterize the mean cluster duration of the Markovian Hawkes process by simplifying the underlying probabilistic combinatorics that have been recently uncovered}. Specifically, we prove that the mean duration is equal to a certain expectation of a $\mathsf{Borel}(\rho)$ random variable (Theorem~\ref{meanIdThm1}), and, by consequence, we obtain an exact expression for the mean duration  (Corollary~\ref{corIdThm1}).

\begin{theorem}\label{meanIdThm1}
Let $B \sim \mathsf{Borel}(\rho)$. Then, the mean duration of the Markovian Hawkes cluster is such that
\begin{align}
\E{\tau}
&=
\left(1-\rho\right)
\E{B \cdot \mathsf{H}_B}
-
1
,
\end{align}
where $\mathsf{H}_j = \sum_{i=1}^j 1/i$ is the $j$th Harmonic number for each $j \in \mathbb{Z}_+$.
\end{theorem}

\begin{corollary} \label{corIdThm1}
The mean duration of the Markovian Hawkes cluster is given by
\begin{align}
\E{\tau}
&=
\int_0^\rho
\frac{e^{-u} - e^{-\rho}}{\rho e^{-\rho} - u e^{-u}}\mathrm{d}u
.
\end{align}
\end{corollary}

\noindent The remainder of this note is devoted to proving this pair. 
In \S\ref{modelSec}, we recast the model via probabilistic combinatorics; in \S\ref{lemmaSec}, we establish associated technical lemmas; in \S\ref{proofSec}, we prove our main results accordingly.

\section{Model Definitions}\label{modelSec}






\subsection{Background Combinatorial Objects}


We begin with a detour. Let us first introduce a pair of (deterministic) combinatorial objects, starting with a well-known extension of parking functions \citep[see, e.g.,][]{stanley2002polytope,yan2015parking}.

\begin{definition}[Generalized Parking Function]\label{gpfDef1}
For $m, n \in \mathbb{Z}_+$ with $m \leq n$, $\pi \in \mathbb{Z}_+^m$ is a \emph{generalized parking function} of length $m$ over $n$ preferences if and only if it is such that, when sorted such that $\pi_{(1)} \leq \dots \leq \pi_{(m)}$, each $\pi_{(i)} \leq n - m + i$ for $i \in \{1, \dots, m\}$.\hfill \Halmos 
\end{definition}

The intuition for generalized parking functions is exactly as the name suggests: for any $m \leq n$, the collection of $\pi \in \mathbb{Z}_+^m$ satisfying Definition~\ref{gpfDef1} contains all the preferences by $m$ drivers over $n$ parking spaces such that, if each driver $i \in \{1, \dots, m\}$ successively proceeds from space 1 to $n$ and takes either space $\pi_i$ or the first available after it, all cars are able to park.
A ``classic'' parking function is obtained when $m = n$. 



In either generalized or classic form, a parking function is inherently ordered by the sequence in which cars arrive, but, as evidenced by the sorting-based condition in Definition~\ref{gpfDef1}, the ability of all drivers to find parking is the same under any ordering of the preferences. This equivalence suggests a possible reduction of information simply to the count of preferences, which we give now for the classic $m = n$ setting.

\begin{definition}[Preference Counts; Cumulative Preference Counts] \label{def:prefs}
For $n \in \mathbb{Z}_+$, $\kappa \in \mathbb{N}^n$ is a vector of \emph{preference counts} if and only if $\sum_{i=1}^n \kappa_i = n$ and $\sum_{j=1}^i \kappa_j \ge i$ for $i \in \{1, \dots, n\}$. Likewise, $\kappa^+ \in \mathbb{Z}^n$ is a vector of \emph{cumulative preference counts} if and only if $1 \leq \kappa_1^+ \leq \dots \leq \kappa_n^+ \leq n$ and $\kappa^+_i \geq i$ for $i \in \{1, \dots, n\}$, where the bijection $\kappa_i^+ = \sum_{j=1}^i \kappa_j$ for $i \in \{1,\dots,n\}$ connects the two concepts.
\hfill \Halmos 
\end{definition}


Definition~\ref{def:prefs} can be readily generalized to the $m < n$ setting, but the  $m=n$ case is sufficient for our remaining pursuit. These preference counts do indeed match the setting of equally enumerated cars and spaces, as the vector defined component-wise by $\kappa_i(\pi) = |\{\ell : \pi_\ell = i \}|$ is a length-$n$ preference count for any length-$n$ parking function $\pi$. Connections beyond parking functions abound: the cumulative preference counts can be quickly recognized as equivalent to Dyck paths on the upper-diagonal half of the $n \times n$ lattice, which are Catalan-numbered. Hence, there are $C_n = (2n)!/(n+1)!n!$ length-$n$ preference counts, by comparison to the $(n+1)^{n-1}$ length-$n$ parking functions.
The non-cumulative preference counts can similarly be thought of in Catalan context as encoding the descendant structure in planar trees, which is a connection we will leverage further in the Hawkes cluster context. 




 



\subsection{Probabilistic Combinatorics Perspectives on the Hawkes Cluster}



Recently, parking functions have been studied with a  probabilistic lens, where these combinatorial objects are sampled uniformly from the set of those with the same length \citep[see, e.g.,][]{diaconis2017probabilizing,kenyon2023parking}. Even more recently, these same objects have been seen to be hidden within Hawkes clusters. \citet{daw2023conditional} first proved that, given the cluster size, the compensator-transformed cluster points are conditionally uniform on a certain polytope that can be parsimoniously partitioned by uniformly random parking functions, and \citet{daw2026asymmetries} then found that uniformly random parking functions also appear in the cluster's original temporal space through the same conditioning on the cluster size.  
In this \emph{end-state conditioning}, the Hawkes cluster model is constructed in opposite order from time: first the total size, then the tree of which points respond to which, and then finally the time between points. In Definition~\ref{def:combx-duration}, we give this construction in its simplest form here for this study of the Hawkes process with its simplest kernel.


\begin{definition}[Combinatorial Cluster Duration] \label{def:combx-duration}
Let $\tau \geq 0$ be constructed as follows:
\begin{enumerate}[i)]
    \item Sample $N \sim \mathrm{Borel}(\rho)$ as the \emph{size} of the cluster. If $N=1$, then $\tau = 0$.
    \item Given $N = n + 1$ for $n \in \mathbb{Z}_+$, sample a length-$n$ preference count $\kappa$ according to the distribution 
    \begin{align}
\PP{\kappa_1 = k_1, \dots, \kappa_n = k_n}
&=
\frac{1}{(n+1)^{n-1}}{n \choose k_1, \dots, k_n}
=
\frac{n!}{(n+1)^{n-1}\prod_{i=1}^n k_i !}
.
\label{kappaDist}
\end{align}
    \item Given $\kappa$ with cumulative form $\kappa^+$, sample $Z_i \sim \mathrm{Exp}(\kappa^+_i -i  + 1)$ independently for each $i \in \{1, \dots, n\}$.
\end{enumerate}
Then, let $\tau := \sum_{i=1}^{N-1} Z_i$ be the \emph{combinatorially constructed cluster duration}.   
\hfill \Halmos
\end{definition}

By Theorem 1 of \citet{daw2026asymmetries} (or Theorem 6 therein for the most specific case), the duration $\tau$ constructed in Definition~\ref{def:combx-duration} is equivalent to the $\tau$ defined by the intensity-based model construction from \eqref{intensityDefEq}. Accordingly, we will base from this construction of $\tau$ in the proofs of our main results. By comparison to \citet{daw2023conditional} or \citet{daw2026asymmetries}, Definition~\ref{def:combx-duration} skips over the dependence on a random parking function $\pi$ and proceeds straight to a random preference count $\kappa$, rather than having the latter being determined by the former. In doing so, we focus on the interpretability of the model via the tree given by $\kappa$, and the distribution of $\kappa$ in \eqref{kappaDist} follows directly from $\kappa_i(\pi) = |\{\ell : \pi_\ell = i \}|$ for a uniformly random $\pi$.



\section{Two Technical Lemmas for Random Preference Counts}\label{lemmaSec}

To leverage this $\kappa$-first construction in Definition~\ref{def:combx-duration}, let us develop two distributional results for the random preference count. The first is a distribution of the next component of $\kappa$ given the preceding components.

\begin{lemma}\label{kappaCPMF1}
For $n \in \mathbb{Z}_+$, let $\kappa$ be a length-$n$ preference count with distribution given by~\eqref{kappaDist}. Then, the conditional distribution of $\kappa_{m+1}$ given $\kappa_1$ through $\kappa_m$ has probability mass function (PMF) given by
\begin{align}
\PP{\kappa_{m+1} = k_{m+1} \mid \kappa_1 = k_1, \dots, \kappa_m = k_m}
&=
\frac{k_{m+1}^+ - m}{k_m^+ - m + 1}
{n - k_m^+ \choose k_{m+1}}
\frac{\left(n-m\right)^{n-k_{m+1}^+-1}}{(n-m+1)^{n-k_m^+-1}}
,
\end{align}
for each $m < n$, meaning that $\kappa_{m+1} | (\kappa_1, \dots, \kappa_m) \stackrel{\mathsf{D}}{=} X_m + Y_m$, where, given $\kappa_m^+$,
\begin{align}
X_m | Y_m
\sim 
\mathsf{Bin}\left(
n - \kappa_m^+ - Y_m
,
\frac{1}{n-m+1}
\right)
\qquad
\text{ and }
\qquad
Y_m \sim 
\mathsf{Bern}\left(\frac{n-\kappa_m^+}{(n-m)(\kappa_m^+ - m + 1)} \right)
,
\end{align}
with $(\kappa_1 - 1) \sim \mathsf{Bin}\left(n-1, 1/(n+1)\right)$.
\end{lemma}
\proof{Proof.}
To start, let us first find the probability of the event that $\kappa_1 = k_1, \kappa_2 = k_2, \dots, \kappa_m = k_m$. 
Consider the possibilities for the underlying parking function $\pi$ which gives $\kappa = \kappa(\pi)$.
Thinking first of $\kappa_1 = k_1$, there are ${n \choose k_1}$ ways to arrange these $k_1$ occurrences of the value $1$ within $\pi$. Moving on, this leaves ${n - k_1 \choose k_2}$ ways to have $\kappa_2 = k_2$, ${n - k_1 - k_2 \choose k_3}$ ways to have $\kappa_3 = k_3$, and so on, for $\prod_{i=1}^m {n - k_{i-1}^+ \choose k_i}$ ways to place the corresponding amounts of the values 1 through $m$ in the vector of length $n$, where we take $k_0^+ = 0$. Then, there remain $n - k_m^+$ unassigned coordinates within the length-$n$ parking function on which the $n-m$ remaining values must be placed. Hence, this remainder constitutes a generalized parking function of length $n - k_m^+$ over $n-m$ preferences. By, e.g., Theorem 2.3 of \citet{kenyon2023parking}, there are $(k_m^+-m+1)(n-m+1)^{n - k_m^+-1}$ such generalized parking functions. Hence, because $\pi$ is uniformly random and because there are $(n+1)^{n-1}$ parking functions of length $n$, we thus have
\begin{align}
\PP{\kappa_1 = k_1, \kappa_2 = k_2, \dots, \kappa_m = k_m}
&=
(k_m^+-m+1)
\frac{(n-m+1)^{n - k_m^+-1}}{(n+1)^{n-1}}
\prod_{i=1}^m {n - k_{i-1}^+ \choose k_i}
.
\end{align}
Then, by conditioning and simplifying, we now find
\begin{align}
\PP{\kappa_{m+1} = k_{m+1} \mid \kappa_1 = k_1, \dots, \kappa_m = k_m}
&=
\frac{
(k_{m+1}^+-m)
(n-m)^{n - k_{m+1}^+-1}
}
{
(k_m^+-m+1)
(n-m+1)^{n - k_m^+-1}
}
{n - k_{m}^+ \choose k_{m+1}}
,
\end{align}
which achieves the stated conditional PMF of $\kappa_{m+1}$ given $\kappa_1$ through $\kappa_m$.
Finally, one can quickly verify that $X_m + Y_m$ satisfies the same conditional PMF given $\kappa_m^+$, and this check completes the proof.
\hfill\Halmos\endproof

Now, for our second technical lemma, let us use the conditional distribution from Lemma~\ref{kappaCPMF1} to find the expected value of $1/(\kappa_{m+1}^+ - m)$ as a recurrence relation indexed by $m$.


\begin{lemma}\label{meanRecur1}
For $n \in \mathbb{Z}_+$, let $\kappa$ be a length-$n$ preference count with distribution given by~\eqref{kappaDist}. Then, for each $m \in \{0, 1, \dots, n - 1\}$,
\begin{align}
\E{\frac{1}{\kappa_{m+1}^+ - m} }
&=
\left(1 + \frac{1}{n-m}\right)\E{\frac{1}{\kappa_m^+ - m + 1}}
-
{n \choose m}\frac{(m+1)^{m-1}}{(n+1)^{n-1}}(n-m)^{n-m-1}
,
\end{align}
where $\PP{\kappa_0^+ = 0} = 1$.
\end{lemma}
\proof{Proof.} Let us start by first proving the following conditional statement:
\begin{align}
\E{\frac{1}{\kappa_{m+1}^+  - m} ~\bigg|~ \kappa_1, \dots, \kappa_m}
&=
\frac{1}{\kappa_m^+ - m + 1}\left(1 + \frac{1}{n-m}\right)
-
\left(\frac{n-m}{n-m+1}\right)^{n-m-1}\mathbf{1}\left\{ \kappa_m^+ = m\right\}
\label{condMeanEq1}
.
\end{align}
From the conditional distribution of $\kappa_{m+1}$ in Lemma~\ref{kappaCPMF1}, we can recognize that $\kappa_{m+1}$ only has a positive probability of being equal to 0 if $\kappa_m^+ > m$; otherwise, if $\kappa_m^+ = m$, then $\kappa_{m+1} \geq 1$. Furthermore, by definition, we have that $\kappa_{m+1} \leq n - \kappa_m^+$. Hence, we can write the conditional expectation of $1/(\kappa_{m+1} + \kappa_m^+ - m)$ as
\begin{align}
\sum_{k=0}^{n-\kappa_m^+}
\frac{1}{\kappa_m^+ - m + 1}
{n - \kappa_m^+ \choose k}
\frac{\left(n-m\right)^{n-k-\kappa_{m}^+-1}}{(n-m+1)^{n-k_m^+-1}}
\nonumber
-
\frac{\left(n-m\right)^{n-m-1}}{(n-m+1)^{n-m-1}}
\mathsf{1}\{\kappa_m^+ = m\}
.
\end{align}
Seeing as the shift term is precisely the lattermost quantity in the stated result, let us now focus on evaluating the sum. Moving forward terms which do not depend on $k$, we can apply the binomial theorem and find
\begin{align*}
\frac{1}{(n-m)(\kappa_m^+ - m + 1)(n-m+1)^{n-\kappa_m^+-1}}
\sum_{k=0}^{n-\kappa_m^+}
{n - \kappa_m^+ \choose k}
\left(n-m\right)^{n-\kappa_{m}^+-k}
&=
\frac{n-m+1}{(n-m)(\kappa_m^+ - m + 1)}
,
\end{align*}
and, together with the shift term, this now simplifies to the stated expression in \eqref{condMeanEq1}.

We can now apply total probability to complete the proof of the original statement. Naturally, \eqref{condMeanEq1} yields
\begin{align}
\E{\frac{1}{\kappa_{m+1}^+ - m} }
&=
\left(1 + \frac{1}{n-m}\right)\E{\frac{1}{\kappa_m^+ - m + 1}}
-
\left(\frac{n-m}{n-m+1}\right)^{n-m-1}\PP{\kappa_m^+ = m}
,
\end{align}
and thus all that is left to evaluate is $\PP{\kappa_m^+ = m}$. For this, we can simply observe that, for $\kappa_m^+$ to be equal to $m$, then the underlying parking function can actually be split into two different smaller parking functions: one of length $m$, which contains the values 1 through $m$, and one of length $n - m$, which contains the values $m+1$ through $n$. Because there are $(m+1)^{m-1}$ parking functions of length $m$, $(n-m+1)^{n-m-1}$ parking functions of length $n-m$, and ${n \choose m}$ ways to arrange the two  smaller parking functions into one parking function of length $n$, we have that
$
\PP{\kappa_m^+ = m}
=
{n \choose m}\frac{(m+1)^{m-1}}{(n+1)^{n-1}}(n-m+1)^{n-m-1}
$,
and, with this expression, now we can simplify to the stated result.
\hfill\Halmos\endproof

The recurrence relation of the expectations in Lemma~\ref{meanRecur1} now immediately yields the following solution.

\begin{corollary}\label{meanRecurSolve1}
For $m < n$, let $\kappa$ be a length-$n$ preference count with distribution given by~\eqref{kappaDist}. Then,
\begin{align}
\E{\frac{1}{\kappa_{m+1}^+ - m} }
&=
\frac{n+1}{n-m}
-
\frac{1}{n-m}
\sum_{i=0}^m
{n \choose i}
\left(\frac{i+1}{n+1}\right)^{i-1}
\left(\frac{n-i}{n+1}\right)^{n-i}
.
\end{align}
\end{corollary}



\section{Proofs of Main Results}\label{proofSec}

With this knowledge of $\kappa$ in hand, we can now prove our main results, Theorem \ref{meanIdThm1} and Corollary \ref{corIdThm1}. 

\subsection{Proof of Theorem \ref{meanIdThm1}}

\proof{Proof.}Let us start the analysis of the mean of $\tau$ through end-state conditioning on the cluster's size:
\begin{align}
\E{\tau}
=
\sum_{n=0}^\infty \E{\tau \mid N = n+1}\PP{N=n+1}
=
\sum_{n=1}^\infty \sum_{m=0}^{n-1} \E{\frac{1}{\kappa_{m+1}^+-m}}  (n+1)^{n-1} \frac{\rho^n e^{-\rho(n+1)}}{n!} 
,
\label{thisOne}
\end{align}
where one should notice that we have dropped the $n=0$ index from the summation as $\E{\tau \mid N = 1} = 0$ by definition. We can now see that Corollary~\ref{meanRecurSolve1} contains exactly the means of the exponential random variables in step (iii) of Definition~\ref{def:combx-duration}, which are the only remaining expectations in \eqref{thisOne}. By applying Corollary~\ref{meanRecurSolve1} and canceling like terms, we can now split this expression into two separate sums,
\begin{align}
\E{\tau}
&=
\sum_{n=1}^\infty 
\sum_{m=0}^{n-1}
\left(
\frac{n+1}{n-m}
-
\frac{1}{n-m}
\sum_{i=0}^m
{n \choose i}
\left(\frac{i+1}{n+1}\right)^{i-1}
\left(\frac{n-i}{n+1}\right)^{n-i}
\right)
(n+1)^{n-1} \frac{\rho^n e^{-\rho(n+1)}}{n!}
\\
&=
\sum_{n=1}^\infty 
\mathsf{H}_n (n+1)^{n} \frac{\rho^n e^{-\rho(n+1)}}{n!}
-
\sum_{n=1}^\infty
\sum_{m=0}^{n-1}
\sum_{i=0}^m
\frac{\left({i+1}\right)^{i-1}
\left({n-i}\right)^{n-i}}{n-m}
 \frac{\rho^n e^{-\rho(n+1)}}{i!(n-i)!}
 .
 \label{twoSumEq1}
\end{align}
Let us analyze each sum separately. 

Starting with the latter, we can first quickly recognize that by exchanging the order of the summations indexed by $m$ and by $i$, we obtain another Harmonic number, $\mathsf{H}_{n-i}$. Then, we next exchange the order of summations over $i$ and over $n$, and, after doing so, we can extract terms that depend on $i$ alone so that the innermost summation is universally a function of the difference between $n$ and $i$. Together, this yields
\begin{align*}
\sum_{n=1}^\infty
\sum_{m=0}^{n-1}
\sum_{i=0}^m
\frac{\left({i+1}\right)^{i-1}
\left({n-i}\right)^{n-i}}{n-m}
 \frac{\rho^n e^{-\rho(n+1)}}{i!(n-i)!}
 \nonumber
 &=
 \sum_{n=1}^\infty
\sum_{i=0}^{n-1}
\mathsf{H}_{n-i}
\left({i+1}\right)^{i-1}
\left({n-i}\right)^{n-i}
 \frac{\rho^n e^{-\rho(n+1)}}{i!(n-i)!}
 \\
  &=
 \sum_{i=0}^\infty
 \frac{\left(\rho({i+1})\right)^{i} }{(i+1)!}
 \rho e^{-\rho(i+1)}
\sum_{n=i+1}^{\infty}
\mathsf{H}_{n-i}
 \frac{\left(\rho(n-i)\right)^{n-i} e^{-\rho(n-i)}}{\rho(n-i)!}
 ,
\end{align*}
and at this point we can recognize that we have two autonomous sums. By re-indexing the right summation to be in terms of $\ell = n - i$, we find that it is equivalent to
\begin{align}
\sum_{\ell=1}^{\infty}
\mathsf{H}_{\ell}
 \frac{\left(\rho\ell\right)^{\ell} e^{-\rho\ell}}{\rho \ell!}
 &=
 \sum_{\ell=1}^{\infty}
\ell \mathsf{H}_{\ell}
 \frac{\left(\rho\ell\right)^{\ell-1} e^{-\rho\ell}}{ \ell!}
 =
 \E{B \mathsf{H}_B}
 ,
\end{align}
for $B \sim \mathsf{Borel}(\rho)$. Then, the left summation can similarly be observed to be 
\begin{align}
 \sum_{i=0}^\infty
 \frac{\left(\rho({i+1})\right)^{i} }{(i+1)!}
 \rho e^{-\rho(i+1)}
 &=
 \rho
  \sum_{j=1}^\infty
 \frac{\left(\rho j\right)^{j-1} }{j!}
 e^{-\rho j}
 =
 \rho
 .
\end{align}
Hence, the latter sum in \eqref{twoSumEq1} simplifies to $\rho \E{B \mathsf{H}_B}$.

Now, for the former sum in \eqref{twoSumEq1}, we can again identify an expectation of a function of $B$ after some quick manipulation. By multiplying and dividing by $n+1$, we can recognize that this is the truncated expectation of $B \mathsf{H}_{B-1}$ when $B$ is restricted to being at least 2:
\begin{align}
\sum_{n=1}^\infty 
\mathsf{H}_n (n+1)^{n} \frac{\rho^n e^{-\rho(n+1)}}{n!}
&=
\sum_{n=1}^\infty 
(n+1)\mathsf{H}_n  \frac{\left(\rho(n+1)\right)^n e^{-\rho(n+1)}}{(n+1)!}
=
\E{B \mathsf{H}_{B-1} \mathbf{1}\{B \geq 2\}}
.
\end{align}
Hence, \eqref{twoSumEq1} now simplifies to
\begin{align}
\E{\tau}
&=
\sum_{n=1}^\infty 
\mathsf{H}_n (n+1)^{n} \frac{\rho^n e^{-\rho(n+1)}}{n!}
-
\sum_{n=1}^\infty
\sum_{m=0}^{n-1}
\sum_{i=0}^m
\frac{\left({i+1}\right)^{i-1}
\left({n-i}\right)^{n-i}}{n-m}
 \frac{\rho^n e^{-\rho(n+1)}}{i!(n-i)!}
\\
&=
\E{B \mathsf{H}_{B-1} \mathbf{1}\{B \geq 2\}}
-
\rho \E{B \mathsf{H}_{B}}
.
\end{align}
Finally, let us now recognize that, at $B = 1$, then
$
B \mathsf{H}_{B-1} \mathbf{1}\{B \geq 2\}
-
\rho B \mathsf{H}_{B}
=
-
\rho
$,
whereas for any $B \geq 2$, we instead find that
$
B \mathsf{H}_{B-1} \mathbf{1}\{B \geq 2\}
-
\rho B \mathsf{H}_{B}
=
B (\mathsf{H}_{B-1} - \mathsf{H}_B)
+
(1-\rho) B \mathsf{H}_{B}
=
-1
+
(1-\rho) B \mathsf{H}_{B}
$.
Thus, we have by linearity of expectation that
$\E{B \mathsf{H}_{B-1} \mathbf{1}\{B \geq 2\}
-
\rho B \mathsf{H}_{B}}
=
(1-\rho)\E{B \mathsf{H}_\mathsf{B}}
-
1$.
\hfill\Halmos\endproof

\subsection{Proof of Corollary \ref{corIdThm1}}


\proof{Proof.} First, by the Eulerian representation of the $j$th Harmonic number as $\mathsf{H}_j = \int_0^1 (1-x^j)/(1-x) \mathrm{d}x$ \citep[e.g.,][sec. 2.2]{lagarias2013euler}, Theorem~\ref{meanIdThm1} can be re-expressed as
\begin{align}
\E{\tau}
&=
(1-\rho)\E{B \mathsf{H}_B} - 1
=
(1-\rho)\int_0^1\frac{\E{B} - \E{B x^B}}{1-x}\mathrm{d}x
-
1
.
\end{align}
By leveraging the fact that the PGF of the Borel distribution can be expressed through the principal branch of the Lambert-W function ($\mathsf{W}_0(\cdot)$) \citep[e.g.,][eq. 3.1]{corless1996lambert}, we can write
\begin{align}
\E{B x^B}
&=
x \frac{\partial}{\partial x}
\left(
\E{x^B}
\right)
=
x \frac{\partial}{\partial x}
\left(
-\frac{1}{\rho}\mathsf{W}_0\left(-x \rho e^{-\rho}\right)
\right)
=
-
\frac{1}{\rho}
\left(
\frac{\mathsf{W}_0\left(-x \rho e^{-\rho}\right)}
{1 + \mathsf{W}_0\left(-x \rho e^{-\rho}\right)}
\right)
;
\end{align}
hence the integral may be equivalently written as
\begin{align}
\int_0^1\frac{\E{B} - \E{B x^B}}{1-x}\mathrm{d}x
&=
\int_0^1\frac{1}{1-x}\left(\frac{1}{1-\rho} +  \frac{1}{\rho}
\left(
\frac{\mathsf{W}_0\left(-x \rho e^{-\rho}\right)}
{1 + \mathsf{W}_0\left(-x \rho e^{-\rho}\right)}
\right)\right)\mathrm{d}x
.
\end{align}

To further simplify this integral, let us apply a change of variable. Let $u = - \mathsf{W}_0(-x\rho e^{-\rho})$. Then, at $x=0$, $u = -\mathsf{W}_0(0) = 0$, and at $x=1$, $u = -\mathsf{W}_0(-\rho e^{-\rho}) = \rho$. Moreover, by the defining property of the Lambert-W, $\mathsf{W}_0(z) e^{\mathsf{W}_0(z)} = z$, we can express $x$ in terms of $u$ via $-u e^{-u} = \mathsf{W}_0(-x\rho e^{-\rho}) e^{\mathsf{W}_0(-x\rho e^{-\rho})} = - x \rho e^{-\rho}$ implying that $x = u e^{-u} / (\rho e^{-\rho})$. By consequence, we also have that $\mathrm{d}x = \mathrm{d}u (e^{-u} - u e^{-u})/(\rho e^{-\rho})$. Thus, the integral over $x$ is equivalent to
\begin{align}
\int_0^1\frac{1}{1-x}\left(\frac{1}{1-\rho} +  \frac{1}{\rho}
\left(
\frac{\mathsf{W}_0\left(-x \rho e^{-\rho}\right)}
{1 + \mathsf{W}_0\left(-x \rho e^{-\rho}\right)}
\right)\right)\mathrm{d}x
&=
\int_0^\rho
\frac{e^{-u} - u e^{-u}}{\rho e^{-\rho} - u e^{-u}}
\left(
\frac{1}{1-\rho}
-
\frac{u}{\rho(1-u)}
\right)
\mathrm{d}u
\\
&=
\frac{1}{\rho(1-\rho)}
\int_0^\rho
\frac{\rho e^{-u} - u e^{-u}}{\rho e^{-\rho} - u e^{-u}} \mathrm{d}u
.
\end{align}
Returning to the expression for the mean duration, we now have
\begin{align}
\E{\tau}
&=
\frac{1}{\rho}
\int_0^\rho
\frac{\rho e^{-u} - u e^{-u}}{\rho e^{-\rho} - u e^{-u}} 
\mathrm{d}u
-
1
=
\frac{1}{\rho}
\int_0^\rho
\frac{\rho e^{-u} - \rho e^{-\rho}}{\rho e^{-\rho} - u e^{-u}} 
\mathrm{d}u
,
\end{align}
and this immediately simplifies to the stated result.
\hfill\Halmos\endproof

\section{Discussion}\label{concSec}

Methodologically, let us note that there may be other proof techniques that could produce these results: for instance, the \citet{moller2005perfect} ODE for $F$ could possibly be used to obtain Corollary~\ref{corIdThm1}, and the $\mathsf{Borel}(\rho)$ identity in Theorem~\ref{meanIdThm1} could then likely be reverse-engineered from the proof of Corollary~\ref{corIdThm1} here. However, we believe that one advantage of the presented end-state conditioning approach is that it preserves interpretability for as long as possible, and we hope that this perspective could potentially provide a way forward for other kernels or other performance metrics. One hurdle for non-exponential kernels will be the lack of memorylessness, which is what allows step (iii) of Definition~\ref{def:combx-duration} to simply depend on the cumulative preference counts. Nevertheless, even in this special and simple setting, Theorem~\ref{meanIdThm1} and Corollary~\ref{corIdThm1} can provide insight. For instance, using either result, one can show that $\E{\tau}/2\log\E{N} \to 1$ as $\rho \to 1$, which demonstrates how much more extreme the growth of the expected cluster size is relative to the growth of its expected duration, and $\E{\tau} \sim -2\log(1-\rho)$ may be a useful asymptotic understanding on its own.







\section*{Acknowledgments}


{We gratefully acknowledge the support of this research from the National Science Foundation (NSF) though grants CMMI-2441387 (AD and AG) and DGE-1842487 (AG).}

\section*{Disclosure of AI Usage}

Generative AI has been used for coding support (such as to enumerate the combinatorial objects while developing conjectures and running simulation experiments to validate results), for reference searches (such as finding a proper citation for the Eulerian Harmonic identity), and for a pre-submission review. 
Notably, the AI-based pre-submission review included a publicly available Lean formalization via PaperFormalizer \citep{tan2026paperformalizer}.
All exposition, statements, and proofs are solely by the authors. Moreover, the authors take responsibility for the correctness of all contents of the paper.

\bibliographystyle{informs2014} 
\bibliography{Bibliography}

\end{document}